\documentclass[a4paper,fleqn]{cas-sc}

\def\tsc#1{\csdef{#1}{\textsc{\lowercase{#1}}\xspace}}
\tsc{WGM}
\tsc{QE}
\tsc{EP}
\tsc{PMS}
\tsc{BEC}
\tsc{DE}
\usepackage[utf8]{inputenc}
\usepackage[english]{babel}
\usepackage{float}
\usepackage[numbers]{natbib}
\begin{document}
\shorttitle{}
\shortauthors{}

\title [mode = title]{The inverse of the (alternating) infinite sum of the reciprocal of the weighted sum for generalized Fibonacci sub-sequences}                      
\tnotemark[1]

\tnotetext[1]{This paper is supported by National Natural Science Foundation of China(Grant No. 12471494) and Natural Science Foundation of Sichuan Province (2024NSFSC2051). The corresponding author is Professor Qunying Liao.}

\author[]{Yongkang Wan}[style=chinese]

\author[]{Zhonghao {Liang}}[style=chinese]

\author[]{Qunying {Liao}$^*$}[style=chinese]
\nonumnote{$^*$Y.Wan, Z.Liang and Q.Liao are with the College of Mathematical Science, Sichuan Normal University, Chengdu 610066, China(email:2475636261@qq.com;liangzhongh0807@163.com;qunyingliao@sicnu.edu.cn)}

\begin{abstract}
In this paper, for the generalized Fibonacci sequence $\left\{W_n\left(a,b,p,q\right)\right\}$, by using elementary methods and techniques, we give the asymptotic estimation values of $\left(\sum\limits_{k=n}^{\infty}\frac{1}{\sum\limits_{i=0}^{t}s_{i}W_{mk+l_i}}\right)^{-1}$ and $\left(\sum\limits_{k=n}^{\infty}\frac{\left(-1\right)^k}{\sum\limits_{i=0}^{t}s_{i}W_{mk+l_i}}\right)^{-1}$, respectively. In particular,  for some special $a,b,p,q,m,t,s_i$ and $l_i\left(0\leq i\leq t \right)$, Theorem \ref{theorem 3.1} is just Theorems 2.1, 2.5-2.6 in \cite{A22} given by Yuan et al.
\end{abstract}




\maketitle

\section{Introduction}
At the beginning of the 12th century, Fibonacci \cite{A1} proposed the famous $\mathrm{Fibonacci}$ sequence $\left\{F_n\right\}$, which is derived from a linear recurrence relation
	$$F_0=0,F_1=1,F_n=F_{n-1}+F_{n-2}\left(n\geq 2\right),$$
	and the corresponding Binet formula is defined as
	$$F_n=\frac{\alpha^n-\beta^n}{\alpha-\beta}\left(n\geq 0\right),$$
	where $\alpha=\frac{1+\sqrt{5}}{2},\beta=\frac{1-\sqrt{5}}{2}$. In 1965, Horadam \cite{A2} proposed a new sequence by  preserving the recurrence relation and altering the first two terms of the $\mathrm{Fibonacci}$ sequence, i.e.,  the generalized $\mathrm{Fibonacci}$ sequence $\left\{W_n\left(a,b,p,q\right)\right\}\left(a,b,p,q\in \mathbb{Z}\right)$, which is derived from the linear recurrence relation
		$$W_0=a,W_1=b,W_n=pW_{n-1}+qW_{n-2}\left(n\geq 2\right).$$
		The corresponding Binet formula is defined as
		$$W_n=c_1\alpha^n-c_2\beta^n\left(n\geq 0\right),$$
		where $c_1=\frac{b-a\beta}{\alpha-\beta},c_2=\frac{b-a\alpha}{\alpha-\beta},\alpha=\frac{p+\sqrt{p^2+4q}}{2}$ and $\beta=\frac{p-\sqrt{p^2+4q}}{2}$.

Since the $\mathrm{Fibonacci}$ sequence has good properties, it's important in theory and  applications, such as in combinatorics [3-4] and cryptography [5-6]. Up to now, the research for the various properties of the (generalized) $\mathrm {Fibonacci}$ sequence has attracted much attention [7-22]. Specifically, in recent years, many scholars considered the asymptotic estimation problem of the reciprocal infinite sum of (generalized) Fibonacci sequences (the definition of asymptotic estimation see Definition \ref{definition 2.1}) and obtained some results [18-22]. For example, in 2021, for the Fibonacci sequence $\left\{F_n\right\}$, Lee \cite {A20} et al. gave the following two asymptotic estimation values,
	$$\left(\sum\limits_{k=n}^{\infty}\frac{1}{F_{k}}\right)^{-1}\sim F_{n-2},$$
	and
$$\left(\sum\limits_{k=n}^{\infty}\frac{1}{F_{mk-l}}\right)^{-1}\sim F_{mn-l}-F_{m\left(n-1\right)-l},$$
where $m,l$ are positive integers  with $l\leq m-1 $.
In 2025, for a special class of  generalized $\mathrm{Fibonacci} $ sequences $\left\{W_n\left(0,1,A,B\right)\right\}$
$\left(A\in \mathbb{Z}^+,B=\pm 1\right)$, Yuan et al.\cite{A22} obtained the following three asymptotic estimation results, 
$$\left(\sum\limits_{k=n}^{\infty}\frac{1}{W_{mk}\left(0,1,A,B\right)}\right)^{-1}\sim W_{mn}\left(0,1,A,B\right)-W_{m\left(n-1\right)}\left(0,1,A,B\right)
,$$
$$\begin{aligned}
	\left(\sum\limits_{k=n}^{\infty}\frac{1}{W_{mk}\left(0,1,A,B\right)+W_{mk+l}\left(0,1,A,B\right)}\right)^{-1}\sim& W_{mn+l}\left(0,1,A,B\right)-W_{m\left(n-1\right)+l}\left(0,1,A,B\right)\\
	&+W_{mn}\left(0,1,A,B\right)-W_{m\left(n-1\right)}\left(0,1,A,B\right)
	,
\end{aligned}$$ 
and 
$$\begin{aligned}
	\left(\sum\limits_{k=n}^{\infty}\frac{1}{\sum\limits_{i=0}^{l}W_{mk+i}\left(0,1,A,B\right)}\right)^{-1}\sim & \frac{1}{\alpha-1}\left(W_{mn+l+1}\left(0,1,A,B\right)-W_{m\left(n-1\right)+l+1}\left(0,1,A,B\right)\right.\\
	&-W_{mn}\left(0,1,A,B\right)+\left.W_{m\left(n-1\right)}\left(0,1,A,B\right)
	\right),
\end{aligned}$$
where $m$ and $l $ are positive integers.

Motivated by the  above works, we 
consider the most general form of the generalized Fibonacci sequence $\left\{W_n\left(a,b,p,q\right)\right\}$ and obtain two types of asymptotic estimation results which generalize  the asymptotic estimation results of Yuan et al. \cite{A22} in 2025. 

For the convenience,  we denote $W_n$ by $W_n\left(a,b,p,q\right)$.

This paper is organized as follows. In Section 2, we give the definitions of the asymptotic estimation and the Big O notation, as well as some necessary lemmas. In Sections 3-4, we obtain the asymptotic estimation values of $\left(\sum\limits_{k=n}^{\infty}\frac{1}{\sum\limits_{i=0}^{t}s_{i}W_{mk+l_i}}\right)^{-1}$ and $\left(\sum\limits_{k=n}^{\infty}\frac{\left(-1\right)^k}{\sum\limits_{i=0}^{t}s_{i}W_{mk+l_i}}\right)^{-1}$, respectively. In Section 5, we conclude the whole paper.
\section{Preliminaries}
In this section, we present the definitions of the asymptotic estimation and the Big O notation, as well as some necessary lemmas. Firstly, we provide the definition of the asymptotic estimation of convergent series as follows.
\newdefinition{definition}{Definition}
\begin{definition}[\cite{A10}]\label{definition 2.1}
	For the convergent series $\sum\limits_{k=1}^{\infty} A_k$, if there exists a function $B_n $ such that
	$$\lim\limits_{n\rightarrow \infty}\left(\left(\sum\limits_{k=n}^{\infty} A_k\right)^{-1}-B_n\right)=0,$$
	then $B_n$ is referred to as the asymptotic estimation of $\left(\sum\limits_{k=n}^{\infty} A_k\right)^{-1}$, and denoted by $\left(\sum\limits_{k=n}^{\infty} A_k\right)^{-1}\sim B_n$.
\end{definition}

Secondly, we review the definition and some properties of the Big O notation.

\begin{definition}[\cite{A15}]\label{definition 2.2}
	Let $f\left(x\right)$ and $g\left(x\right)\in \mathbb{R}[x]$ with
	$\lim\limits_{x\rightarrow x_0}f\left(x\right)=0$ and $\lim\limits_{x\rightarrow x_0}g\left(x\right)=0$, respectively, where $x_0\in \mathbb{R}$ 
	. If $x$ lies in a certain punctured neighborhood of $x_0 $, and there exists a constant $A>0$ such that
	$$\left|\frac{f\left(x\right)}{g\left(x\right)}\right|\leq A,$$
	then  $\frac{f\left(x\right)}{g\left(x\right)}$ is called an bounded quantity as $x \rightarrow x_0$, and denoted by $f\left(x\right)=O\left(g\left(x\right)\right)\left(x\rightarrow x_0\right).$
\end{definition}

By Definition 2, it's easy to prove  the following

\newtheorem{lemma}{Lemma}\label{lemma 2.1}
\begin{lemma}
	If $x$ is a real number  with $|x|<1 $, then the following statements are true,
	
	$\left(1 \right)$ for any constant $C$ and positive integer $n$, 
	$$O\left(Cx^n\right)=O\left(x^n\right) \ \text{and}\ C\cdot O\left(x^n\right)=O\left(x^n\right);$$
	
	$\left(2\right)$ for any positive integers $n$ and $m$ with $n \leq m $, 
	$O\left(x^n\right)+O\left(x^m\right)=O\left(x^n\right);$
	
	$\left(3\right)$ for any real number $y$ with $|y|>1 $, positive integers $n$ and $m$ , 
	$O\left(\frac{x^n}{y^n}\right)\cdot \frac{1}{y^m}=O\left(\frac{x^n}{y^{n+m}}\right).$
\end{lemma}

By the proof of Theorem 1 in \cite{A16}, we have the following 

\begin{lemma}[\cite{A16}]\label{lemma 2.2}
	If $x$ is a real number with $| x |<1 $, then we have
	$$\frac{1}{1\pm x}=1+O\left(x\right).$$
\end{lemma}

Finally, to prove our main results, we need the following
\begin{lemma}\label{lemma 2.3}
	For any positive integers $m, k, t, p$, integers $q,l_i$ and natural numbers $s_i\left(i=0,1,\ldots,t\right)$, where $l_i\geq 1-m$, and $\left(s_0,s_1,\cdots,s_t\right)\in \mathbb{N}^{t+1}\backslash\left\{\mathbf{0}\right\}$, we have
	$$
	\frac{1}{\sum\limits_{i=0}^{t}s_{i}W_{mk+l_i}}=\frac{\frac{1}{\alpha^{mk}}+O\left(\frac{\beta^{mk}}{\alpha^{2mk}}\right)}{c_1\sum\limits_{i=0}^{t}s_{i} \alpha^{l_i}}.
	$$
\end{lemma}
\newproof{pf}{Proof}
\begin{pf}
	From the Binet formula of the generalized $\mathrm{Fibonacci}$ sequence $\left\{W_n\left(a,b,p,q\right)\right\}$, it's easy to know that
	$$
	\begin{aligned}
		\frac{1}{\sum\limits_{i=0}^{t}s_{i}W_{mk+l_i}}&=\frac{1}{\sum\limits_{i=0}^{t}s_{i}\left(c_1\alpha^{mk+l_i}-c_2\beta^{mk+l_i}\right)}\\
		&=\frac{1}{c_1\sum\limits_{i=0}^{t}s_{i} \alpha^{mk+l_i}-c_2\sum\limits_{i=0}^{t}s_{i}\beta^{mk+l_i}}\\
		&=\frac{1}{c_1\sum\limits_{i=0}^{t}s_{i} \alpha^{mk+l_i}}\left(1-\frac{c_2}{c_1}\left(\frac{\beta}{\alpha}\right)^{mk}\frac{\sum\limits_{i=0}^{t}s_{i}\beta^{l_i}}{\sum\limits_{i=0}^{t}s_{i}\alpha^{l_i}}\right)^{-1}.
	\end{aligned}
	$$
By $p>0 $, it's easy to prove that $\left|\frac{\beta}{\alpha}\right|<1$. Therefore, for any positive integer $k$, we have
$$\lim\limits_{k\rightarrow \infty}\left(\frac{\beta}{\alpha}\right)^{k}=0,$$
thus by the definition, it's easy to know that there exists a sufficiently large positive integer $k$ such that
	$$\left|\frac{c_2}{c_1}\left(\frac{\beta}{\alpha}\right)^{mk}\frac{\sum\limits_{i=0}^{t}s_{i}\beta^{l_i}}{\sum\limits_{i=0}^{t}s_{i}\alpha^{l_i}}\right|<1.$$
	Furthermore, by Lemmas 1-2, we have
	$$
	\begin{aligned}
		\frac{1}{\sum\limits_{i=0}^{t}s_{i}W_{mk+l_i}}&=\frac{1}{c_1\sum\limits_{i=0}^{t}s_{i} \alpha^{mk+l_i}}\left(1+O\left(\frac{\beta^{mk}}{\alpha^{mk}}\right)\right)\\
		&=\frac{\frac{1}{\alpha^{mk}}+O\left(\frac{\beta^{mk}}{\alpha^{2mk}}\right)}{c_1\sum\limits_{i=0}^{t}s_{i} \alpha^{l_i}}
		.
	\end{aligned}
	$$
	
	And then we complete the proof of Lemma \ref{lemma 2.3}.
\end{pf}

\section{The asymptotic estimation value of $\left(\sum\limits_{k=n}^{\infty}\frac{1}{\sum\limits_{i=0}^{t}s_{i}W_{mk+l_i}}\right)^{-1}$}
In this section, for the generalized Fibonacci sequence $\left\{W_n\left(a,b,p,q\right)\right\}$, we give the asymptotic estimation value of $\left(\sum\limits_{k=n}^{\infty}\frac{1}{\sum\limits_{i=0}^{t}s_{i}W_{mk+l_i}}\right)^{-1}$, i.e., we prove the following 

\newtheorem{theorem}{Theorem}
\begin{theorem}\label{theorem 3.1}
For any positive integers $m, k, t, p$, integers $q,l_i$ and natural numbers $s_i\left(i=0,1,\ldots,t\right)$, where $p^2+2q-2<p\sqrt{p^2+4q},l_i\geq 1-m$ and $\left(s_0,s_1,\cdots,s_t\right)\in \mathbb{N}^{t+1}\backslash\left\{\mathbf{0}\right\}$, we have
$$\left(\sum_{k=n}^{\infty}\frac{1}{\sum\limits_{i=0}^{t}s_{i}W_{mk+l_i}}\right)^{-1}\sim \sum\limits_{i=0}^{t}s_{i}\left(W_{mn+l_i}-W_{m\left(n-1\right)+l_i}\right)  .$$
\end{theorem}
\begin{pf}
By Lemma \ref{lemma 2.3}, we have
$$
\sum_{k=n}^{\infty}\frac{1}{\sum\limits_{i=0}^{t}s_{i}W_{mk+l_i}}=\frac{1}{c_1\sum\limits_{i=0}^{t}s_{i} \alpha^{l_i}}\left(\sum_{k=n}^{\infty}\frac{1}{\alpha^{mk}}+\sum_{k=n}^{\infty}O\left(\frac{\beta^{mk}}{\alpha^{2mk}}\right)\right).
$$
Note that
$$
\sum_{k=n}^{\infty}\frac{1}{\alpha^{mk}}=\sum_{i=0}^{\infty}\frac{1}{\alpha^{m\left(n+i\right)}}=
\frac{1}{\alpha^{mn}}\sum_{i=0}^{\infty}\left(\frac{1}{\alpha^{m}}\right)^{i}=\frac{1}{\alpha^{mn}}\cdot\frac{1}{1-\frac{1}{\alpha^{m}}}\\
=\frac{\alpha^{m}}{\alpha^{mn}\left(\alpha^{m}-1\right)},
$$
and then by Lemma \ref{lemma 2.2}, it's easy to get
$$\sum_{k=n}^{\infty}O\left(\frac{\beta^{mk}}{\alpha^{2mk}}\right)=O\left(\frac{\beta^{mn}}{\alpha^{2mn}}\right),$$
thus
$$\begin{aligned}
	\left(\sum_{k=n}^{\infty}\frac{1}{\sum\limits_{i=0}^{t}s_{i}W_{mk+l_i}}\right)^{-1}&=\left(\frac{1}{c_1\sum\limits_{i=0}^{t}s_{i} \alpha^{l_i}}\left(\frac{\alpha^{m}}{\alpha^{mn}\left(\alpha^{m}-1\right)}+O\left(\frac{\beta^{mn}}{\alpha^{2mn}}\right)\right)\right)^{-1}\\
	&=\left(\frac{1}{c_1\sum\limits_{i=0}^{t}s_{i} \alpha^{l_i}}\cdot\frac{\alpha^{m}}{\alpha^{mn}\left(\alpha^{m}-1\right) }\left(1+O\left(\frac{\beta^{mn}}{\alpha^{mn}}\right)\right)\right)^{-1}\\
	&=c_1\sum\limits_{i=0}^{t}{s_{i} \alpha^{l_i}}\cdot\frac{\alpha^{mn}\left(\alpha^{m}-1\right) }{\alpha^{m}}\left(1+O\left(\frac{\beta^{mn}}{\alpha^{mn}}\right)\right)^{-1}.
\end{aligned}$$
Now by Definition \ref{definition 2.2}, there exists a constant $A>0$  such that
$$\left|O\left(\frac{\beta^{mn}}{\alpha^{mn}}\right)\right|<A\left|\left(\frac{\beta}{\alpha}\right)^{mn}\right|.$$
Note that $p>0$, it's easy to prove that $\left|\frac{\beta}{\alpha}\right|<1$, and so, for any positive integer $n$, we have
$$\lim\limits_{n\rightarrow \infty}A\left(\frac{\beta}{\alpha}\right)^{n}=0,$$
then by the definition, it's easy to know that there exists a sufficiently large positive integer $n$ such that
$$A\left|\left(\frac{\beta}{\alpha}\right)^{mn}\right|<1,$$
i.e.,
$$
	\left|O\left(\frac{\beta^{mn}}{\alpha^{mn}}\right)\right|<1.
$$
Furthermore, by Lemmas 1-2, we have
$$
	\begin{aligned}
		&\left(\sum_{k=n}^{\infty}\frac{1}{\sum\limits_{i=0}^{t}s_{i}W_{mk+l_i}}\right)^{-1}\\
		=&c_1\sum\limits_{i=0}^{t}s_{i} \alpha^{l_i}\cdot\frac{\alpha^{mn}\left(\alpha^{m}-1\right) }{\alpha^{m}}\left(1+O\left(\frac{\beta^{mn}}{\alpha^{mn}}\right)\right)\\
		=&c_1\sum\limits_{i=0}^{t}s_{i} \alpha^{l_i}\left(\alpha^{mn}-\alpha^{m\left(n-1\right)}\right)+c_1\sum\limits_{i=0}^{t}s_{i} \alpha^{l_i}\cdot\frac{\alpha^{m}-1 }{\alpha^{m}}\cdot\alpha^{mn}\cdot  O\left(\frac{\beta^{mn}}{\alpha^{mn}}\right)\\
		=&\sum\limits_{i=0}^{t}s_{i} \left(c_1\alpha^{mn+l_i}-c_1\alpha^{m\left(n-1\right)+l_i}\right) +O\left(\beta^{mn}\right)\\
		=&\sum\limits_{i=0}^{t}s_{i}\left(W_{mn+l_i}+c_2\beta^{mn+l_i}-W_{m\left(n-1\right)+l_i}-c_2\beta^{m\left(n-1\right)+l_i}\right)  +O\left(\beta^{mn}\right)\\
		=&\sum\limits_{i=0}^{t}s_{i}\left(W_{mn+l_i}-W_{m\left(n-1\right)+l_i}\right)  +\sum\limits_{i=0}^{t}s_{i}\left(c_2\beta^{mn+l_i}-c_2\beta^{m\left(n-1\right)+l_i}\right)  +O\left(\beta^{mn}\right).
	\end{aligned}
$$
Note that $p^2+2q-2<p\sqrt{p^2+4q}$, it's easy to prove that $|\beta|<1$, hence,
$$
\begin{aligned}
	&\lim_{n\rightarrow \infty}
	\left(\left(\sum_{k=n}^{\infty}\frac{1}{\sum\limits_{i=0}^{t}s_{i}W_{mk+l_i}}\right)^{-1}-\left(\sum\limits_{i=0}^{t}s_{i}\left(W_{mn+l_i}-W_{m\left(n-1\right)+l_i}\right)\right)\right)\\
	=&\lim_{n\rightarrow \infty}
	\left(\sum\limits_{i=0}^{t}s_{i}\left(c_2\beta^{mn+l_i}-c_2\beta^{m\left(n-1\right)+l_i}\right) +O\left(\beta^{mn}\right)\right)\\
	=&\lim_{n\rightarrow \infty}
	\left(c_2\beta^{m\left(n-1\right)}\sum\limits_{i=0}^{t}s_{i}\left( \beta^{m+l_i}- \beta^{l_i}\right)+O\left(\beta^{mn}\right)\right)\\
	=&0.
\end{aligned}
$$
Therefore, by Definition \ref{definition 2.1}, we get
$$\left(\sum_{k=n}^{\infty}\frac{1}{\sum\limits_{i=0}^{t}s_{i}W_{mk+l_i}}\right)^{-1}\sim \sum\limits_{i=0}^{t}s_{i}\left(W_{mn+l_i}-W_{m\left(n-1\right)+l_i}\right) .$$

This completes the proof of Theorem \ref{theorem 3.1}.
\end{pf}

In particular, for some special $t, l_i$ and $s_i$, we can immediately get the following Corollaries 1-3.

\newtheorem{corollary}{Corollary}
\begin{corollary}\label{corollary 3.1}
For any positive integers $ m,k,p,s_0$ and integers $q,l_0$, where $p^2+2q-2<p\sqrt{p^2+4q}$ and $l_0\geq 1-m$, we have
$$\left(\sum_{k=n}^{\infty}\frac{1}{s_{0}W_{mk+l_0}}\right)^{-1}\sim s_{0}  W_{mn+l_0}-s_{0} W_{m\left(n-1\right)+l_0}.$$
\end{corollary}

\newdefinition{rmk}{Remark}
\begin{rmk}
	$\left(1\right)$
	By taking $a=l_0=0$ and $b=p=q=s_0=m=1$ in Corollary \ref{corollary 3.1}, the corresponding  result is just Theorem 1.1 in \cite{A20};
	
	$\left(2\right)$
	By taking $a=0,b=p=q=s_0=1$ and $1-m\leq l_0<0$ in Corollary \ref{corollary 3.1}, the corresponding  result is just Theorem 1.3 in \cite{A20};
	
	$\left(3\right)$
	By taking $a=l_0=0,b=p=q=s_0=m=1$ in Corollary \ref{corollary 3.1}, the corresponding  result is just Theorem 2.1 in \cite{A22}.
\end{rmk}

\begin{corollary}\label{corollary 3.2}
For any positive integers $ m,k,p$, integers $q,l_i$ and natural numbers $s_i\left(i=0,1\right)$, where $p^2+2q-2<p\sqrt{p^2+4q},l_i\geq 1-m$ and $\left(s_0,s_1\right)\in \mathbb{N}^{2}\backslash\left\{\mathbf{0}\right\}$, we have
$$\left(\sum_{k=n}^{\infty}\frac{1}{s_{0}W_{mk+l_0}+s_{1}W_{mk+l_1}}\right)^{-1}\sim s_{0}  W_{mn+l_0}-s_{0} W_{m\left(n-1\right)+l_0}+s_{1}  W_{mn+l_1}-s_{1} W_{m\left(n-1\right)+l_1}.$$
\end{corollary}

\begin{rmk}
	By taking $a=0,b=1,p=A\in \mathbb{Z}^+,q=\pm 1,s_0=s_1=1$ and $l_0=0,l_1=d\in \mathbb{Z}^+$ in Corollary \ref{corollary 3.2}, the corresponding  result is just Theorem 2.5 in \cite{A22}.
\end{rmk}

\begin{corollary}\label{corollary 3.3}
	For any positive integers $ m,k,p,t$, integers $q$ and $i$, where $p^2+2q-2<p\sqrt{p^2+4q}$, we have
	$$\left(\sum_{k=n}^{\infty}\frac{1}{\sum\limits_{i=0}^{t}W_{mk+i}}\right)^{-1}\sim \frac{1}{\alpha-1}\left( W_{mn+t+1}- W_{mn}- W_{m\left(n-1\right)+t+1}+  W_{m\left(n-1\right)}\right).$$
\end{corollary}
\begin{pf}
 From the proof of Theorem \ref{theorem 3.1}, we know that
$$\begin{aligned}
	&\left(\sum_{k=n}^{\infty}\frac{1}{\sum\limits_{i=0}^{t}W_{mk+i}}\right)^{-1}\\
	=&\sum\limits_{i=0}^{t} c_1\alpha^{mn+i}-\sum\limits_{i=0}^{t} c_1\alpha^{m\left(n-1\right)+i}+O\left(\beta^{mn}\right)\\
	=&c_1\alpha^{mn}\cdot\sum_{i=0}^{t}\alpha^{i}-c_1\alpha^{m\left(n-1\right)}\cdot \sum_{i=0}^{t}\alpha^{i}+O\left(\beta^{mn}\right)\\
	=&c_1\alpha^{mn}\cdot \frac{\alpha^{t+1}-1}{\alpha-1}-c_1\alpha^{m\left(n-1\right)}\cdot \frac{\alpha^{t+1}-1}{\alpha-1}+O\left(\beta^{mn}\right)\\
	=&\frac{1}{\alpha-1}\left(c_1\alpha^{mn+t+1}-c_1\alpha^{mn}-c_1\alpha^{m\left(n-1\right)+t+1}+c_1\alpha^{m\left(n-1\right)}\right)+O\left(\beta^{mn}\right).
\end{aligned}$$
Note that for any positive integers $e$ and $l$, we have
$$W_{me+l}=c_1\alpha^{me+l}-c_2\beta^{me+l},$$
and so
$$\begin{aligned}\left(\sum_{k=n}^{\infty}\frac{1}{\sum\limits_{i=0}^{t}W_{mk+i}}\right)^{-1}=&\frac{1}{\alpha-1}\left(W_{mn+t+1}-W_{mn}-W_{m\left(n-1\right)+t+1}+W_{m\left(n-1\right)}\right)\\&+\frac{1}{\alpha-1}\left(c_2\beta^{mn+t+1}-c_2\beta^{mn}-c_2\beta^{m\left(n-1\right)+t+1}+c_2\beta^{m\left(n-1\right)}\right)+O\left(\beta^{mn}\right).
\end{aligned}$$
Note that $p^2+2q-2<p\sqrt{p^2+4q}$, it's easy to prove that $|\beta|<1$, hence,
$$
\begin{aligned}
	&\lim_{n\rightarrow \infty}
	\left(\left(\sum_{k=n}^{\infty}\frac{1}{\sum\limits_{i=0}^{t}W_{mk+i}}\right)^{-1}-\left(\frac{1}{\alpha-1}\left(W_{mn+t+1}-W_{mn}-W_{m\left(n-1\right)+t+1}+W_{m\left(n-1\right)}\right)\right)\right)\\
	=&\lim_{n\rightarrow \infty}
	\left(\frac{1}{\alpha-1}\left(c_2\beta^{mn+t+1}-c_2\beta^{mn}-c_2\beta^{m\left(n-1\right)+t+1}+c_2\beta^{m\left(n-1\right)}\right)+O\left(\beta^{mn}\right)\right)\\
	=&\lim_{n\rightarrow \infty}
	\left(\frac{c_2}{\alpha-1}\cdot \beta^{m\left(n-1\right)}\left(\beta^{m+t+1}-\beta^m-\beta^{t+1}+1\right)+O\left(\beta^{mn}\right)\right)\\
	=&0.
\end{aligned}
$$
Therefore, by Definition \ref{definition 2.1}, we know that
$$\left(\sum_{k=n}^{\infty}\frac{1}{\sum\limits_{i=0}^{t}W_{mk+i}}\right)^{-1}\sim \frac{1}{\alpha-1}\left(W_{mn+t+1}-W_{mn}-W_{m\left(n-1\right)+t+1}+W_{m\left(n-1\right)}\right).$$

And then we complete the proof of Corollary \ref{corollary 3.3}.
\end{pf}

\begin{rmk}
	By taking $a=0,b=1,p=A\in \mathbb{Z}^+,q=\pm 1$ and $t=l$ in Corollary \ref{corollary 3.3}, the corresponding  result is just Theorem 2.6 in \cite{A22}.
\end{rmk}

\section{The asymptotic estimation value of $\left(\sum\limits_{k=n}^{\infty}\frac{\left(-1\right)^k}{\sum\limits_{i=0}^{t}s_{i}W_{mk+l_i}}\right)^{-1}$}
In this section, for the generalized Fibonacci sequence $\left\{W_n\left(a,b,p,q\right)\right\}$, we give the asymptotic estimation value of $\left(\sum\limits_{k=n}^{\infty}\frac{\left(-1\right)^k}{\sum\limits_{i=0}^{t}s_{i}W_{mk+l_i}}\right)^{-1}$, i.e., we prove the following 


\begin{theorem}\label{theorem4.1}
For any positive integers $m, k, t, p$, integers $q,l_i$ and natural numbers $s_i\left(i=0,1,\ldots,t\right)$, where $p^2+2q-2<p\sqrt{p^2+4q},l_i\geq 1-m$ and $\left(s_0,s_1,\cdots,s_t\right)\in \mathbb{N}^{t+1}\backslash\left\{\mathbf{0}\right\}$, we have
$$\left(\sum_{k=n}^{\infty}\frac{\left(-1\right)^k}{\sum\limits_{i=0}^{t}s_{i}W_{mk+l_i}}\right)^{-1}\sim \left(-1\right)^n\sum\limits_{i=0}^{t}s_{i} \left( W_{mn+l_i}+ W_{m\left(n-1\right)+l_i}\right).$$
\end{theorem}
\begin{pf}
By Lemma \ref{lemma 2.3}, we have
$$
\sum_{k=n}^{\infty}\frac{\left(-1\right)^k}{\sum\limits_{i=0}^{t}s_{i}W_{mk+l_i}}=\frac{1}{c_1\sum\limits_{i=0}^{t}s_{i} \alpha^{l_i}}\left(\sum_{k=n}^{\infty}\frac{\left(-1\right)^k}{\alpha^{mk}}+\sum_{k=n}^{\infty}\left(-1\right)^k O\left(\frac{\beta^{mk}}{\alpha^{2mk}}\right)\right).
$$
Note that
$$
\sum_{k=n}^{\infty}\frac{\left(-1\right)^k}{\alpha^{mk}}=\sum_{i=0}^{\infty}\frac{\left(-1\right)^{n+i}}{\alpha^{m\left(n+i\right)}}=
\frac{\left(-1\right)^n}{\alpha^{mn}}\sum_{i=0}^{\infty}\left(\frac{-1}{\alpha^m}\right)^{i}=\frac{\left(-1\right)^n}{\alpha^{mn}}\cdot\frac{1}{1-\frac{-1}{\alpha^{m}}}\\
=\frac{\left(-1\right)^{n}\alpha^{m}}{\alpha^{mn}\left(\alpha^{m}+1\right)},
$$
and then by Lemma 2, we can get
$$\sum_{k=n}^{\infty}\left(-1\right)^kO\left(\frac{\beta^{mk}}{\alpha^{2mk}}\right)=\sum_{k=n}^{\infty}O\left(\frac{\beta^{mk}}{\alpha^{2mk}}\right)=O\left(\frac{\beta^{mn}}{\alpha^{2mn}}\right),$$
thus,
$$\begin{aligned}
	\left(\sum_{k=n}^{\infty}\frac{\left(-1\right)^k}{\sum\limits_{i=0}^{t}s_{i}W_{mk+l_i}}\right)^{-1}&=\left(\frac{1}{c_1\sum\limits_{i=0}^{t}s_{i} \alpha^{l_i}}\left(\frac{\left(-1\right)^n\alpha^{m}}{\alpha^{mn}\left(\alpha^{m}+1\right)}+O\left(\frac{\beta^{mn}}{\alpha^{2mn}}\right)\right)\right)^{-1}\\
	&=\left(\frac{1}{c_1\sum\limits_{i=0}^{t}s_{i} \alpha^{l_i}}\cdot\frac{\left(-1\right)^n\alpha^{m}}{\alpha^{mn}\left(\alpha^{m}+1\right) }\left(1+O\left(\frac{\beta^{mn}}{\alpha^{mn}}\right)\right)\right)^{-1}\\
	&=\left(-1\right)^{n} c_1\sum\limits_{i=0}^{t}{s_{i} \alpha^{l_i}}\cdot\frac{\alpha^{mn}\left(\alpha^{m}+1\right) }{\alpha^{m}}\left(1+O\left(\frac{\beta^{mn}}{\alpha^{mn}}\right)\right)^{-1}.
\end{aligned}$$
Now from the proof of Theorem \ref{theorem 3.1}, we can get
$$\left|O\left(\frac{\beta^{mn}}{\alpha^{mn}}\right)\right|<1,$$
and then by Lemma \ref{lemma 2.1},
$$
	\begin{aligned}
		&\left(\sum_{k=n}^{\infty}\frac{\left(-1\right)^k}{\sum\limits_{i=0}^{t}s_{i}W_{mk+l_i}}\right)^{-1}\\
		=&\left(-1\right)^{n} c_1\sum\limits_{i=0}^{t}{s_{i} \alpha^{l_i}}\cdot\frac{\alpha^{mn}\left(\alpha^{m}+1\right) }{\alpha^{m}}\left(1+O\left(\frac{\beta^{mn}}{\alpha^{mn}}\right)\right)\\
		=&\left(-1\right)^{n}c_1\sum\limits_{i=0}^{t}s_{i} \alpha^{l_i}\left(\alpha^{mn}+\alpha^{m\left(n-1\right)}\right)+\left(-1\right)^{n}c_1\sum\limits_{i=0}^{t}s_{i} \alpha^{l_i}\cdot\frac{\alpha^{m}+1 }{\alpha^{m}}\cdot\alpha^{mn}\cdot  O\left(\frac{\beta^{mn}}{\alpha^{mn}}\right)\\
		=&\left(-1\right)^n\sum\limits_{i=0}^{t}s_{i} \left( c_1\alpha^{mn+l_i}+ c_1\alpha^{m\left(n-1\right)+l_i}\right)+O\left(\beta^{mn}\right)\\
		=&\left(-1\right)^n\sum\limits_{i=0}^{t}s_{i} \left( W_{mn+l_i}+c_2\beta^{mn+l_i}+ W_{m\left(n-1\right)+l_i}+ c_2\beta^{m\left(n-1\right)+l_i}\right)+O\left(\beta^{mn}\right)\\
		=&\left(-1\right)^n\sum\limits_{i=0}^{t}s_{i} \left( W_{mn+l_i}+ W_{m\left(n-1\right)+l_i}\right)+\left(-1\right)^n\sum\limits_{i=0}^{t}s_{i} \left( c_2\beta^{mn+l_i}+ c_2\beta^{m\left(n-1\right)+l_i}\right)+O\left(\beta^{mn}\right).
	\end{aligned}
$$
Note that $p^2+2q-2<p\sqrt{p^2+4q}$, it's easy to prove that $|\beta|<1$, hence,
$$
\begin{aligned}
	&\lim_{n\rightarrow \infty}
	\left(\left(\sum_{k=n}^{\infty}\frac{\left(-1\right)^k}{\sum\limits_{i=0}^{t}s_{i}W_{mk+l_i}}\right)^{-1}-\left(-1\right)^n\sum\limits_{i=0}^{t}s_{i}\left(  W_{mn+l_i}+ W_{m\left(n-1\right)+l_i}\right)\right)\\
	=&\lim_{n\rightarrow \infty}
	\left(\left(-1\right)^n\sum\limits_{i=0}^{t}s_{i} \left( c_2\beta^{mn+l_i}+ c_2\beta^{m\left(n-1\right)+l_i}\right)+O\left(\beta^{mn}\right)\right)\\
	=&\lim_{n\rightarrow \infty}
	\left(c_2\beta^{m\left(n-1\right)}\left(\left(-1\right)^n\sum\limits_{i=0}^{t}s_{i} \beta^{m+l_i}+\left(-1\right)^n\sum\limits_{i=0}^{t}s_{i} \beta^{l_i}\right)+O\left(\beta^{mn}\right)\right)\\
	=&0.
\end{aligned}
$$
Therefore, by Definition \ref{definition 2.1}, we can get
$$\left(\sum_{k=n}^{\infty}\frac{\left(-1\right)^k}{\sum\limits_{i=0}^{t}s_{i}W_{mk+l_i}}\right)^{-1}\sim \left(-1\right)^n\sum\limits_{i=0}^{t}s_{i} \left( W_{mn+l_i}+ W_{m\left(n-1\right)+l_i}\right).$$

From the above, we complete the proof of Theorem \ref{theorem4.1}.
\end{pf}

In particular, for some special $t, l_i$ and $s_i$, we can immediately get the following Corollaries 4-6.

\begin{corollary}\label{corollary 4.1}
For any positive integers $ m,k,p,s_0$, and integers $q,l_0$, where $p^2+2q-2<p\sqrt{p^2+4q}$ and $l_0\geq 1-m$, we have
$$\left(\sum_{k=n}^{\infty}\frac{\left(-1\right)^k}{s_{0}W_{mk+l_0}}\right)^{-1}\sim \left(-1\right)^ns_{0}\left(  W_{mn+l_0}+ W_{m\left(n-1\right)+l_0}\right).$$
\end{corollary}

\begin{corollary}\label{corollary 4.2}
For any positive integers $ m,k,p$, integers $q,l_i$ and natural numbers $s_i\left(i=0,1\right)$, where $p^2+2q-2<p\sqrt{p^2+4q},l_i\geq 1-m$ and $\left(s_0,s_1\right)\in \mathbb{N}^{2}\backslash\left\{\mathbf{0}\right\}$, we have
$$\left(\sum_{k=n}^{\infty}\frac{\left(-1\right)^k}{s_{0}W_{mk+l_0}+s_{1}W_{mk+l_1}}\right)^{-1}\sim \left(-1\right)^n\left(s_{0}  W_{mn+l_0}+s_{0} W_{m\left(n-1\right)+l_0}+s_{1}  W_{mn+l_1}+s_{1} W_{m\left(n-1\right)+l_1}\right).$$
\end{corollary}

\begin{corollary}\label{corollary 4.3}
For any positive integers $ m,k,p,t$, integers $q,i$, where $p^2+2q-2<p\sqrt{p^2+4q}$, we have
$$\left(\sum_{k=n}^{\infty}\frac{\left(-1\right)^k}{\sum\limits_{i=0}^{t}W_{mk+i}}\right)^{-1}\sim \frac{\left(-1\right)^{n}}{\alpha-1}\left(W_{mn+t+1}-W_{mn}+W_{m\left(n-1\right)+t+1}-W_{m\left(n-1\right)}\right).$$
\end{corollary}
\begin{pf}
From the proof of Theorem \ref{theorem4.1}, we have
$$\begin{aligned}
	\left(\sum_{k=n}^{\infty}\frac{\left(-1\right)^k}{\sum\limits_{i=0}^{t}W_{mk+i}}\right)^{-1}&=\left(-1\right)^n\left(\sum\limits_{i=0}^{t} c_1\alpha^{mn+i}+\sum\limits_{i=0}^{t} c_1\alpha^{m\left(n-1\right)+i}\right)+O\left(\beta^{mn}\right)\\
	&=\left(-1\right)^{n}c_1\alpha^{mn}\cdot \sum_{i=0}^{t}\alpha^{i}+\left(-1\right)^{n}c_1\alpha^{m\left(n-1\right)}\cdot \sum_{i=0}^{t}\alpha^{i}+O\left(\beta^{mn}\right)\\
	&=\left(-1\right)^{n}c_1\alpha^{mn}\cdot \frac{\alpha^{t+1}-1}{\alpha-1}+\left(-1\right)^{n}c_1\alpha^{m\left(n-1\right)}\cdot \frac{\alpha^{t+1}-1}{\alpha-1}+O\left(\beta^{mn}\right)\\
	&=\frac{\left(-1\right)^{n}}{\alpha-1}\left(c_1\alpha^{mn+t+1}-c_1\alpha^{mn}+c_1\alpha^{m\left(n-1\right)+t+1}-c_1\alpha^{m\left(n-1\right)}\right)+O\left(\beta^{mn}\right).
\end{aligned}$$
Note that for any positive integers $e$ and $l$, we have
$$W_{me+l}=c_1\alpha^{me+l}-c_2\beta^{me+l},$$
and so 
$$\begin{aligned}\left(\sum_{k=n}^{\infty}\frac{\left(-1\right)^k}{\sum\limits_{i=0}^{t}W_{mk+i}}\right)^{-1}=&\frac{\left(-1\right)^{n}}{\alpha-1}\left(W_{mn+t+1}-W_{mn}+W_{m\left(n-1\right)+t+1}-W_{m\left(n-1\right)}\right)\\
	&+\frac{\left(-1\right)^{n}}{\alpha-1}\left(c_2\beta^{mn+t+1}-c_2\beta^{mn}+c_2\beta^{m\left(n-1\right)+t+1}-c_2\beta^{m\left(n-1\right)}\right)+O\left(\beta^{mn}\right).
\end{aligned}$$
Now by $p^2+2q-2<p\sqrt{p^2+4q}$, it's easy to prove that $|\beta|<1$, hence,
$$
\begin{aligned}
	&\lim_{n\rightarrow \infty}
	\left(\left(\sum_{k=n}^{\infty}\frac{\left(-1\right)^k}{\sum\limits_{i=0}^{t}W_{mk+i}}\right)^{-1}-\left(\frac{\left(-1\right)^{n}}{\alpha-1}\left(W_{mn+t+1}-W_{mn}+W_{m\left(n-1\right)+t+1}-W_{m\left(n-1\right)}\right)\right)\right)\\
	=&\lim_{n\rightarrow \infty}
	\left(\frac{\left(-1\right)^{n}}{\alpha-1}\left(c_2\beta^{mn+t+1}-c_2\beta^{mn}+c_2\beta^{m\left(n-1\right)+t+1}-c_2\beta^{m\left(n-1\right)}\right)+O\left(\beta^{mn}\right)\right)\\
	=&\lim_{n\rightarrow \infty}
	\left(\frac{\left(-1\right)^{n}c_2}{\alpha-1}\cdot \beta^{m\left(n-1\right)}\left(\beta^{m+t+1}-\beta^m+\beta^{t+1}-1\right)+O\left(\beta^{mn}\right)\right)\\
	=&0.
\end{aligned}
$$
Therefore, by Definition \ref{definition 2.1}, we have
$$\left(\sum_{k=n}^{\infty}\frac{\left(-1\right)^k}{\sum\limits_{i=0}^{t}W_{mk+i}}\right)^{-1}\sim \frac{\left(-1\right)^{n}}{\alpha-1}\left(W_{mn+t+1}-W_{mn}+W_{m\left(n-1\right)+t+1}-W_{m\left(n-1\right)}\right).$$

From the above, we complete the proof of Corollary \ref{corollary 4.3}.
\end{pf}

\section{Conclusions}
In this paper, we consider the asymptotic estimation for two classes of generalized Fibonacci sub-sequences and obtain the following main results, which generalize  the asymptotic estimation results of Yuan et al. \cite{A22} in 2025. 
\begin{itemize}
	\item The asymptotic estimate of $\left(\sum\limits_{k=n}^{\infty}\frac{1}{\sum\limits_{i=0}^{t}s_{i}W_{mk+l_i}}\right)^{-1}$(Theorem \ref{theorem 3.1}). In particular, we have the following table.
	\begin{center}
		\scriptsize
		\renewcommand{\arraystretch}{1.5}
		\resizebox{\linewidth}{!}{
			\begin{tabular}{|c|c|c|c|c|c|c|c|} 
				\hline
				Our results & $\left(a,b\right) \in \mathbb{Z}^2$  & $\left(p,q\right) \in \mathbb{Z}^2$ & $m \in \mathbb{Z}^+$ & $t \in \mathbb{N}$ & $(s_0,s_1,\ldots,s_t)\in \mathbb{N}^{t}\backslash \left\{\boldsymbol{0}\right\}$ & $(l_0,l_1,\ldots,l_t)\in \mathbb{Z}^{t}$ & Known results \\
				\hline
				
				\multirow{5}{*}{Theorem \ref{theorem 3.1}} & \multirow{5}{*}{$\left(0,1\right)$} & \multirow{2}{*}{$\left(1,1\right)$}  & 1 & 0 & $s_0 = 1$ & $l_0 = 1$ & \cite{A20}Theorem 1.1 \\ 
				\cline{4-8}
				&& & $\mathbb{Z}^+$ & 0 & $s_0 = 1$& $1 - m \leq l_0 < 0$ & \cite{A20}Theorem1.3 \\ 
				\cline{3-8}
				&&\multirow{3}{*}{$\left(\mathbb{Z}^+,\pm1\right)$}& \multirow{3}{*}{$\mathbb{Z}^+$} & 0 & $s_0 = 1$& $l_0 \in \mathbb{Z}^+$ &\cite{A22}Theorem2.1 \\ \cline{5-8}
				&  &&&1& $\left(s_0 ,s_1\right)= \left(1,1\right)$& $\left(l_0,l_1\right) = \left(0,d\right) ,d\in \mathbb{Z}^+$ & \cite{A22}Theorem2.5 \\ 
				\cline{5-8}
				&&&& $\mathbb{N}$ &$(s_0,s_1,\ldots,s_t)=(1,1,\ldots,1)$ & $ (l_0,l_1,\ldots,l_t)= (0, 1, \ldots, t)$ & \cite{A22}Theorem2.6 \\ 
				\hline
			\end{tabular}
		}
	\end{center}
	\item The asymptotic estimate of $\left(\sum\limits_{k=n}^{\infty}\frac{\left(-1\right)^k}{\sum\limits_{i=0}^{t}s_{i}W_{mk+l_i}}\right)^{-1}$(Theorem \ref{theorem4.1}).
\end{itemize}




\begin{thebibliography}{99}  
	\bibitem{A1} Sigler L. Fibonacci's Liber Abaci: A translation into modern English of Leonardo Pisano's book of calculation[M]. Springer New York, 2002.
	\bibitem{A2} Horadam A F. Generating functions for powers of a certain generalised sequence of numbers[J]. Duke Math, 1965, 32: 437–446.
	\bibitem{A3} De Souza J, Curado E M F, Rego-Monteiro M A. Generalized Heisenberg algebras and Fibonacci series[J]. Journal of Physics A: Mathematical and General, 2006, 39(33): 10415-10425.
	\bibitem{A4} Holliday S H, Komatsu T. On the sum of reciprocal generalized Fibonacci numbers[J]. Integers,  2011, 11(4): 441-455.
	\bibitem{A5} Basu M, Prasad B. Long range variations on the Fibonacci universal code[J]. Journal of Number Theory, 2010, 130(9): 1925-1931.
	\bibitem{A6} Atanassov K T, Atanassova L C, Sasselov D D. A new perspective to the generalization of the Fibonacci sequence[J]. The Fibonacci Quarterly, 1985, 23(1): 21-28.
	\bibitem{A7} Zhang W. Some identities involving the Fibonacci numbers[J]. The Fibonacci Quarterly, 1997, 35(3): 225-228.
	\bibitem{A8} Zhang W. On Chebyshev polynomials and Fibonacci numbers[J]. The Fibonacci Quarterly, 2002, 40(5): 424-428.
	\bibitem{A9} Zhang W. Some identities involving the Fibonacci numbers and Lucas numbers[J]. The Fibonacci Quarterly, 2004, 42(2): 149-154.
	\bibitem{A10}Hu J. Some formulate for the Fibonacci and Lucas numbers[J]. Research On Number Theory And Smarandache Nottons, 2010: 119.
	\bibitem{A11}Li Y. On Chebyshev polynomials, Fibonacci polynomials, and their derivatives[J]. Journal of Applied Mathematics, 2014, 2014(1): 451953.
	\bibitem{A12} Liu, J, Zhang, W. Identities for the sums of binomial coefficients of generalized three - periodic Fibonacci sequences[J]. Pure and Applied Mathematics, 2021, 37(1), 38 - 47.
	\bibitem{A13} Ohtsuka H, Nakamura S. On the sum of reciprocal Fibonacci numbers[J]. The Fibonacci Quarterly, 2008, 46(2): 153-159.
	\bibitem{A14} Zhang G J. The infinite sum of reciprocal of the Fibonacci numbers[J]. J. Math. Res. Expo, 2011, 31(6): 1030-1034.
	\bibitem{A15} Wang A Y Z, Zhang F. The reciprocal sums of even and odd terms in the Fibonacci sequence[J]. Journal of Inequalities and Applications, 2015, 2015: 1-13.
	\bibitem{A16} Wu Z, Zhang H. On the reciprocal sums of higher-order sequences[J]. Advances in Difference Equations, 2013, 2013(1): 189-197.
	\bibitem{A17} Yuan P, He Z, Zhou J. On the sum of reciprocal generalized Fibonacci numbers[C]//Abstract and Applied Analysis. Hindawi Limited, 2014.
	\bibitem{A18} Hwang W T, Park J D, Song K. On the reciprocal sum of the fourth power of Fibonacci numbers[J]. Open Mathematics, 2022, 20(1): 1642-1655.
	\bibitem{A19} Lee H H, Park J D. Asymptotic behavior of reciprocal sum of two products of Fibonacci numbers[J]. Journal of Inequalities and Applications, 2020, 2020(1): 91-107.
	\bibitem{A20} Lee H H, Park J D. The limit of reciprocal sum of some subsequential Fibonacci numbers[J]. AIMS Math, 2021, 6(11): 12379-12394.
	\bibitem{A21} Marques D, Trojovský P. The proof of a formula concerning the asymptotic behavior of the reciprocal sum of the square of multiple-angle Fibonacci numbers[J]. Journal of Inequalities and Applications, 2022, 2022(1): 21-37.
	\bibitem{A22} Li H, Yang K, Yuan P. The asymptotic behavior of the reciprocal sum of generalized Fibonacci numbers[J]. Electronic Research Archive, 2025, 33(1): 409-432.
	\bibitem{A23} Chen, J, Yu, C, Jin, L. Mathematical Analysis (3rd ed., Vol. 1)[M]. Beijing: Higher Education Press, 2019.
	\end{thebibliography}



\bio{}

\endbio



\vskip3pc



\end{document}